\documentclass[11pt]{amsart}
\usepackage{amssymb ,amsfonts}
\usepackage{xypic}
\usepackage{amsthm}

\newtheorem{thm}{Theorem}[section]
\newtheorem{cor}[thm]{Corollary}
\newtheorem{pro}[thm]{Proposition}
\newtheorem{lem}[thm]{Lemma}

\theoremstyle{definition}
\newtheorem{defn}[thm]{Definition}

\newtheorem{exmp}[thm]{Example}

\newtheorem{rem}[thm]{Remark}

\newcommand{\mm}{\mathcal{M}}
\newcommand{\nn}{\mathcal{N}}
\newcommand{\inj}{\text{inj} \, }
\newcommand{\proj}{\text{proj} \, }
\newcommand{\im}{\text{im}}
\newcommand{\Ho}{\text{Ho}}
\newcommand{\pt}{\text{pt} ( \mathcal{E} )}
\newcommand{\ee}{\mathcal{E}}

\newcommand{\sss}{\mathcal{S}}
\newcommand{\kk}{\mathcal{K}}
\newcommand{\p}{\text{pro-} }
\newcommand{\g}{\mathcal{K}}
\newcommand{\oo}{\mathcal{O}}
\newcommand{\dd}{\mathcal{D}}
\newcommand{\cc}{\mathcal{C}}

\newcommand{\rarr}{\rightarrow}
\newcommand{\hh}{\mathcal{H}}
\newcommand{\zz}{\mathbb{Z}}
\newcommand{\spec}{\text{spec} \, }
\newcommand{\dfn}{\textbf}
\newcommand{\mdfn}[1]{\dfn{\mathversion{bold}#1}}
\newcommand{\shc}{\mathcal{E}}
\newcommand{\pshc}{\mathcal{P}}
\newcommand{\ch}{\text{ch} }
\newcommand{\colim}{\text{colim}}
\newcommand{\ww}{E}
\newcommand{\hofib}{\text{hofib}}

\makeatletter
\let\c@equation\c@thm
\makeatother \numberwithin{equation}{section}

\title{t-model structures on
chain complexes of presheaves}
\author{H. Fausk}

\thanks{Work done during a visit to the Institut Mittag-Leffler
(Djursholm, Sweden)}

\date{ August 1.~2007}
\address{
Department of Mathematical Sciences NTNU, Trondheim 7034,
 Norway} \email{fausk@math.ntnu.no}

\begin{document}
\begin{abstract}
A tensor model structure is constructed on the category of chain
complexes of presheaves of $R$-modules for a sheaf of rings $R$ in a Grothendieck
topos. If the topos has enough points, then the homotopy category is
equivalent to the derived category. Some t-structures, generalizing the perverse t-structures, are constructed on the derived category.
\end{abstract}

\maketitle

\section{introduction}
There is an injective and a projective model structure on the
category of chain complexes of $R$-modules for a unital ring $R$
\cite[2.3]{hov}. The injective model structure extends to the
category of chain complexes for any Grothendieck abelian category.
This model structure fails to interact well with the tensor
structure; there are too many cofibrant objects. On the other hand,
the projective model structure is well-behaved with respect to the
tensor product, but it does not readily generalize to chain
complexes of sheaves. For certain ringed spaces, Mark Hovey has
constructed a tensor model structure on the category of chain
complexes of sheaves of $R$-modules \cite{hov01}. See also Gillespie \cite{gil}.

The first goal of this paper is to construct a tensor model
structure on the category of chain complexes of \emph{presheaves} of
$R$-modules for any ringed Grothendieck topos. We do this in two
steps. First we generalize the projective model structure to the
category of chain complexes of presheaves of $R$-modules. Then we
localize this model structure with respect to the stalkwise homology
isomorphisms.

The second goal of this paper is to study some t-structures on the
derived category $\dd_R$.
 These
  t-structures ``lift'' to the model
 category of chain complexes of presheaves to form a t-model
 structure
 \cite{tfi}.

We now give a summary. In Section \ref{sec:derived} we
recall the definition of the derived category. In Section
\ref{sec:proj} we generalize the projective model structure to the
category of chain complexes of presheaves of $R$-modules. The weak
equivalences are presheaf homology-isomorphisms. In Section
\ref{sec:tensor} we prove that the projective model structure is
%monoidal.
tensorial.
%These two sections are elementary and
%more efficient proofs can be given using the well documented results
%in the case of an ordinary ring (ringed topos with one morphism).
In Section \ref{sec:quasi} we define quasi-simplicial model
structures and show that the projective model
structure is quasi-simplicial.

In Section \ref{sec:stalkwise} we localize the projective model
category with respect to the stalkwise homology isomorphisms. We call the resulting model
structure the stalkwise model structure. It is a tensor model category. If the ringed topos has
enough points, then the weak equivalences are exactly the sheaf
homology isomorphisms. So under this assumption the homotopy
category of the stalkwise model category is equivalent to the
derived category $\dd_R$.
%Moreover, if the topos has enough points,
%then the stalkwise model structure is Quillen equivalent to both the
%flat and the injective model structures. There are fewer
%cofibrations in the stalkwise model structure than in the other two
%model structures.

In Section
 \ref{sec:t-structures} we construct some families of
 t-structures on $\dd_R$ and show that they
lift to form t-model
 structures on the category of chain complexes.
 The t-structures we consider are
 ``locally'' the standard t-structures. The simplest class
 of t-structures are
 given by assigning to each point in
 the topos a cut-off value in $\zz \cup \{ \pm \infty \}$.
 %Furthermore we refine this by taking the structure of
% the stalk of rings, $R_p$, into account for all all points $p$ of
% the topos.
  In particular, we generalize
 the perverse t-structures on $\dd_R$ and show that they lift to
 t-model structures.
  In Section \ref{sec:ddgeq1} we give
  an explicit description of $(\dd_R)_{\geq 0} $ and
  $(\dd_R)_{\leq 0}$ associated
to particularly well-behaved t-structures.
%We construct the t-model structures on
%the category of chain complexes of presheaves of $R$-modules by
%localizing with respect to suitable classes of maps that are closed
%under nonnegative suspensions. Hence our approach gives an
%alternative to the gluing-techniques \cite{bbd}.

In Section \ref{sec:examples} we consider various examples and
compare our model structures to the flat model structure and to the
injective model structure. In an Appendix we recall, and slightly
extend, Bousfield's cardinality argument. This is needed in Section
\ref{sec:t-structures}.

We assume the reader is familiar with the fundamentals of model
category theory. See for example \cite{hir,hov}.
%The usual projective model structure on the
%category of chain complexes of $R$-modules, for a ring $R$, is
%obtained as a special case of our results (the one morphism site).
%More precisely, we needed a model structure on $\ch (\mm)$ that
%would allow us to form a tensor triangulated homotopy category of
%pro-chain complexes of $R$-modules.

\section{The derived category}
\label{sec:derived} In this section some notation and terminology is
introduced. Let $\cc$ be a (skeletally) small Grothendieck site. Let
$\shc$ be the category of sheaves of sets on $\cc$, and let $\pshc$
be the category of presheaves of sets on $\cc$. Let $ i \colon \shc
\rarr \pshc$ be the inclusion functor. It has a left adjoint, the
sheafification functor. We denote the sheafification functor by
$L^2$ (=$L \circ L$), and the unit of the adjunction by $ \eta
\colon 1 \rarr i \circ L^2$ \cite[II.3.0.5]{sga4.1}. Assume that
$\ee$ has a set of isomorphism classes of points, and let $\pt$
denote this set.

Let $R$ be a sheaf of rings on $\cc$. Let $\mm$ denote the category
of left $R$-modules in $\shc$, and let $ \nn$ denote the category of
left $i R$-modules in $ \pshc$. Both $\mm$ and $\nn$ are abelian
closed tensor categories with units $R$ and $i R$, respectively. The
functor $i$ induces an inclusion functor $i \colon \mm \rarr \nn$,
and the functor $L^2$ induces a left adjoint sheafification functor
$ L^2 \colon \nn \rarr \mm$.

For any object $C$ in $\cc$, let \mdfn{$ R_C$} denote the free
$R$-module in $\nn$ generated by $C$ \cite[IV.11.3.3]{sga4.1}. There
is a natural isomorphism \begin{equation} \label{eq:RC}
 \nn ( R_C , X ) \cong X
(C) .
\end{equation}
%and $ \nn
%( R_C , X ) \cong X (C)$ has a natural module structure over $ \nn (
%R_C , R) \cong R (C)$.
Similarly, $L^2 R_C$ is the free $R$-module in $\mm$ generated by
$C$. Let $\bullet$ be the terminal object in $\cc$. Then $R$ is
isomorphic to $R_{\bullet}$.

\begin{defn}
Let \mdfn{$\ch (\nn ) $} denote the category of chain complexes of
presheaves of $i R$-modules
 on $\cc$, and let \mdfn{$\ch (\mm ) $}
 denote the category of chain complexes of
sheaves of $R$-modules
 on $\cc$. \end{defn} The categories
 {$\ch (\mm ) $} and {$\ch (\nn ) $} are abelian closed tensor
 categories.

\begin{defn}
\label{defn:quasi-iso} A map $f \colon X \rarr Y$ in $\ch (\nn)$ is
a \mdfn{presheaf homology-isomorphism} if
%for each $n \in \zz$, the induced map on homology in $\nn$,
\[ H_n (f) \colon H_n (X) \rarr H_n (Y) \] is an isomorphism,
for each $n \in \zz$.

A map $f \colon X \rarr Y$ in $\ch (\nn)$ is a \mdfn{sheaf
homology-isomorphism} if the sheafification of the induced map on
homology
\[L^2 H_n (f) \colon L^2 H_n (X) \rarr L^2
H_n (Y) \] is an isomorphism, for each $n \in \zz$.
\end{defn}
\begin{defn} \label{defn:stalkwiseiso}
A map $f \colon X \rarr Y$ in $\ch (\nn )$
 is a \mdfn{stalkwise homology-isomorphism} if $ (L^2 f)_p$ is a
homology-isomorphism of chain complexes of $R_p$-modules for all
points $p$ in $\shc$. \end{defn}
%Note that if $ L^2 f$ is a
%homology-isomorphism of presheaves, then $ f$ is a stalkwise
%equivalence. The converse is false.
Let $i$ also denote the inclusion functor $i \colon \ch (\mm) \rarr
\ch (\nn)$. A map $f$ in $\ch (\mm)$ is a presheaf
homology-isomorphisms, sheaf homology-isomorphism, or stalkwise
homology-isomorphism if $i (f)$ has this property.

 If $\ee$ has enough points, then a
map $f$ in $\ch (\nn ) $ is a stalkwise homology-isomorphism if and
only if it is a sheaf homology-isomorphism \cite[IV.6.4.1]{sga4.1}
\begin{defn} \label{defn:derivedcat}
The localization of $ \ch (\mm )$ with respect to the class of all
sheaf homology-isomorphisms is called the \mdfn{derived category} of
chain complexes of sheaves of $R$-modules on $\cc$. It is denoted by
\mdfn{$\dd_R$}.
\end{defn}
The unit  $ X \rarr  i \circ L^2 X$ of the $(L ,i)$-adjunction is a sheaf homology-isomorphism, for all presheaves $X$. So the derived category is equivalent to the localization of $ \ch ( \nn)$ with respect to the class of all sheaf homology-isomorphisms. (The injective model structure, for example,  gives the existence of the derived categories.)

In Section \ref{sec:stalkwise} a tensor model structure on $\ch
(\nn)$ is constructed such that the weak equivalences are the
stalkwise homology-isomorphisms. Its homotopy category is equivalent
to the derived category provided  $\ee$ has enough points.
 We
first describe a preliminary model structure on $\ch (\nn )$ with
weak equivalences the smaller class of presheaf
homology-isomorphisms.

\section{The projective model structure}
\label{sec:proj} The projective model structure on the
category of chain complexes of $R$-modules, for a ring $R$, is generalized  to the
category of chain complexes of presheaves of $R$-modules, for a
sheaf of rings $R$.

We define the cofibrant
  generators.
\begin{defn} \label{defn:IJ}
Let \mdfn{$i_{C,n}$} be the map of chain complexes
\[
\xymatrix{ \cdots \ar[r] \ar[d] & 0 \ar[r] \ar[d] &0 \ar[r] \ar[d] &
R_C \ar[d]^= \ar[r] & 0 \ar[r] \ar[d] & \cdots \ar[d] \\
\cdots \ar[r] & 0 \ar[r] & R_C \ar[r]^= & R_C \ar[r] & 0 \ar[r] &
\cdots }
\] where $C \in \cc$ and the vertical identity map on $R_C$
 is in degree $n$.
 Let \mdfn{$I$} be
the set of all $ i_{C ,n}$ for $ C \in \cc$ and $n \in \zz$.

Let \mdfn{$j_{C,n}$} be the map of chain complexes \[ \xymatrix{
\cdots \ar[r] \ar[d] & 0 \ar[r] \ar[d] & 0 \ar[r] \ar[d] & 0 \ar[d]
\ar[r] & 0 \ar[r] \ar[d] & \cdots \ar[d] \\ \cdots \ar[r] & 0 \ar[r]
& R_C \ar[r]^= & R_C \ar[r] & 0 \ar[r] & \cdots }
\] where $C \in \cc$ and the rightmost copy of $R_C$ is in degree
 $n$. Let
\mdfn{$J$} be the set of all $ j_{C ,n}$ for $ C \in \cc$ and $n \in
\zz$.
\end{defn}
%The sign are introduced to make $I$ and $J$ closed under suspension
%functor.
The following model structure on $\ch (\nn )$ is called the
\mdfn{projective model structure} on the category of chain complexes
of presheaves of $R$-modules. Relative $I$-cell complexes are
defined in \cite[10.5]{hir}.

\begin{thm} \label{thm:projectivemodelstructure}
There is a proper cofibrantly generated model structure on $\ch (\nn
)$. The weak equivalences are presheaf homology-isomorphisms. The
cofibrations are retracts of relative $I$-cell complexes. The
fibrations are levelwise surjective maps of presheaves. The
cofibrant generators, $I$, and the acyclic cofibrant generators,
$J$, are small.
\end{thm}

A map of presheaves $ f \colon X \rarr Y$ is (levelwise) surjective
if $f_n (C)$ is surjective as a map of sets, for all $C \in \cc$ and
$n \in \zz$. This is equivalent to $f$ being epic (in the
categorical sense) in the category $\nn$.
%%What happens for $\mm$???
%The cofibrant objects are levelwise projective, and flat as chain
%complexes.

\begin{proof}
It suffices to check that: (1) $ \inj (I) $ is equal to $\inj(J) \cap
W$, and (2) $\proj ( \inj (J) ) $ is contained in $W$
\cite[11.3.1]{hir}.

We first show that $ \inj J$ is equal to the class of surjective
maps, and that $ \inj I$ is equal to the class of surjective maps
that are also homology-isomorphisms. This gives $(1)$ and the
description of the fibrant objects.

  A map from $j_{C ,n-1}$ to $f \colon X \rarr Y$ is
specified by $ y \in Y_n (C) $, and a lift is given by an element $x
\in X_n (C)$ such that $f_n (C) ( x ) = y $. Hence $f$ has the right
lifting property with respect to $j_{C, n-1}$ if and only if $ f_n
(C)$ is surjective.

We denote the $n$-cycles of $Y$ at $C$, $\ker (Y_n (C) \rarr Y_{n-1}
(C) )$,
 by $ Z(Y (C))_n $. Observe that a map $ f \colon X \rarr
Y$ has the right lifting property with respect to $ i_{C ,n}$ if and
only if the canonical map from $ X_{n +1 } (C) $ to $W_n (C)$ in the
pullback diagram
\[ \xymatrix{ W_n (C) \ar[d] \ar[r] & Y_{n+1} (C) \ar[d] \\
Z ( X (C))_n \ar[r] & Z(Y (C))_n} \] is surjective.

  Assume that $f$ has the right
lifting property with respect to $i_{C ,n-1}$. Then \[ Z(X (C))_n
\rarr Z(Y (C))_n \] is surjective since $0 \times Z(Y (C))_n$ is
contained in $ W_{n-1} (C)$. Hence $H_n (f) (C)$ is surjective and
$W_n (C) \rarr Y_{n+1} (C)$ is surjective. The right lifting
property with respect to both $i_{C , n-1}$ and $i_{C ,n}$ implies
that $f_{n+1} (C) \colon  X_{n+1} (C) \rarr Y_{n+1 }(C) $ is surjective. So the induced
map on boundaries
\[ \im ( X_{n+1}
(C) \rarr X_{n} (C) ) \rarr \im (Y_{n+1} (C) \rarr Y_{n} (C) ) \] is
surjective. Hence $H_n (f) (C)$ is bijective.

We now prove the converse claim. Assume that $f$ is a
homology-isomorphism of presheaves and $f_n (C)$ is
surjective, for all $C \in \cc$ and $n \in \zz$. Given an element $
x \in Z(X (C))_n$ and an element $y \in Y_{n+1} (C)$ such that $d
(y) = f_n (x)$. We need to show that the element $ ( x , y) \in W_n
(C)$ comes from an element in $X_{n +1 } (C)$. By our assumption
there exists an element $x' \in X_{n+1} (C)$ such that $f_{n+1} (C)
( x' ) = y $. Hence $H_n (f)$ sends the cycle $d(x' ) - x $ to zero.
So $d (x' ) - x $ is a boundary since $H_n (f)$ is injective. We
conclude that there exists an element $x'' \in X_{n+1} (C)$ such
that $ d x'' = x$. We get that $f_{n+1} (C) ( x'' ) - y$ is a cycle.
Since $H_{n+1 } (f) (C)$ is surjective there exists an element $
x''' \in X_{n+1 } (C)$ such that $ d x''' = x$ and $f_{n+1} (C) (
x''' ) =y$. This shows that the diagram lifts. So $i_{C, n}$ has the
left lifting property with respect to $f$.

 We now verify (2). Assume that $f
\colon X \rarr Y $ is a map in $\proj ( \inj (J) ) $. Consider the
diagram \[ \xymatrix{ X \ar[d]^f \ar[r]^(0.35){f , \, \text{id}} & Y
\oplus X \ar[d]^{\text{id}_Y
\oplus 0} \\
Y \ar[r]^= & Y . } \] The rightmost vertical map is surjective. So
by our assumption on $f$   the diagram lifts. Hence $H_n (f)$ is injective
for all $n$.

Now consider the diagram
\[
\xymatrix{ X \ar[d]^f \ar[r]^(0.25){\text{id} ,\, 0 } & X \oplus
\text{Tot} ( Y \oplus Y ) \ar[d]^{f \oplus g}
 \\
Y \ar[r]^= & Y } \] where $\text{Tot} (Y \oplus Y) $ is the total
complex associated to the double complex \[ \xymatrix{ \cdots \ar[r]
\ar[d] & Y_{n+1 } \ar[r] \ar[d]^{(-\text{id})^{n+1} }&
 Y_{n} \ar[r] \ar[d]^{(-\text{id})^{n} } & Y_{n-1 } \ar[r]
 \ar[d]^{(-\text{id})^{n-1} } &
 \cdots \ar[d] \\
 \cdots \ar[r] &
Y_{n+1 } \ar[r] &
 Y_{n} \ar[r] & Y_{n-1 } \ar[r] & \cdots .} \]
The map $g \colon \text{Tot} (Y \oplus Y) \rarr Y$ is given by the
identity map on the upper copy of $Y$ in the double complex and by 0
on the lower copy of $Y$. Hence the rightmost map in the diagram is
surjective. By our assumption on $f$ there is a lift in the diagram. Since
the homology of $\text{Tot} ( Y \oplus Y)$ is $0$ we get that $ H_n
(f)$ is surjective for all $n$.

The verification of properness reduces to the category of chain
complexes of $R(C)$-modules for a ring $R(C)$, for each $C \in \cc$.
A tedious verification shows that the pushout of a
homology-isomorphism along a levelwise injective map of chain
complexes is again a homology-isomorphism. A simpler verification
shows that the pullback of a homology-isomorphism along a levelwise
surjective map of chain complexes is again a homology-isomorphism.
Both pushouts and pullbacks in $\ch (\nn)$ are formed levelwise.
Since the cofibrations are levelwise injective, and the fibrations
%,evaluated at any $C\in\cc$,
are levelwise surjective, it follows that the model structure is
proper.

The cofibrant generators are small since evaluation of presheaves at
an object $C$ of $\cc$ commutes with direct sum.
\end{proof}

 We refer to Hovey for an alternative
  description of the cofibrant objects and
the cofibrations \cite[2.3.6, 2.3.8-9]{hov}. Note that the
isomorphism in equation \ref{eq:RC} shows that the presheaf $R_C$ is a
projective object in $\nn$, for each object $C \in \cc$. In fact all
projective objects in $\nn$ are retracts of direct sums of object of
the form $R_C$, where $C \in \cc$.

\begin{rem}
 The projective model structure on $\ch (
\nn )$ is a stable model category. Hence its  homotopy category is a triangulated category by \cite[7.1]{hov}. This
triangulated category is typically different from the derived
category of chain complexes of sheaves of $R$-modules on the site
$\cc$.
\end{rem}

\begin{defn}
The \mdfn{unit interval complex}, denoted \mdfn{$U$}, consists of
two copies of $R$ in degree 0, and one copy of $R$ in degree 1, the
differential is the identity map on the first copy and minus the
identity map on the second copy of $R$.
\end{defn} If $X$ is a cofibrant object, then a cylinder object
for $X$ is given by $X \textstyle\coprod X \rarr X \otimes U \rarr
X$, where $U$ is the unit interval complex. Hence we get the
following.
\begin{lem} If
  $X$ is a cofibrant object and $Y$ an arbitrary object in the
  projective model structure on $\ch ( \nn )$, then $
\Ho ( \ch (\nn ) )( X , Y )$ is isomorphic to the group of chain
homotopy classes of (degree 0) chain maps from $X$ to $Y$. \end{lem}
\begin{proof} This follows since all objects are fibrant in the
 projective model structure on $\ch ( \nn )$. \end{proof}

\begin{rem}
The projective model structure on $\ch ( \nn )$ is certainly well
known. We give an alternative way to construct the projective model
structure on chain complexes of $R$-modules that mesh better with
the language used in other examples of model structures of this
kind.

There is a tensor category of  presheaves on $\cc$ in the category of
chain complexes of abelian groups.  Let
$A$ be a ring object (monoid) in this abelian category. Then the
category of $A$-modules in the category of presheaves of chain
complexes of abelian groups is equivalent to $A \text{-} \ch ( \nn
)$ (see Definition \ref{defn:a-alg}). The projective model structure
on $A \text{-} \ch ( \nn)$, under this identification, is
%on the category of chain complexes of presheaves of chain
%complexes of abelian groups
inherited from the set of right adjoint  functors \[ A \text{-} \ch (
\nn ) \rarr \ch ( \zz ) \] given by evaluating at objects $C$ in $\cc$ and composing  with the
forgetful functor from $A(C)$-modules to $\zz$-modules (at least one
object in each isomorphism class in $\cc$).

A similar model structure on the category of presheaves of
simplicial sets has been given by Benjamin Blander \cite{bla}. See
also Sharon Hollander's model structure on the category of stacks
\cite{hol}.

A more elaborate class of examples of this type are the strict model
structures on diagram spectra. These examples are most naturally
studied in an enriched setting.
%Let $\cc$ be a Grothendick site enriched in a
%model category $\ttt$, and let $\mm$ be a model category ``model
%theoretically enriched'' over $\ttt$. Then the category of
%presheaves (preserving the enrichment) with values in $\mm$
%is ``model theoretically enriched'' over $\ttt$.
See \cite{mj}, \cite{ren} and \cite{ss}.

%The ``sheaf
%condition'' is given by require $ \mm ( C, X ) \rarr \mm ( F, X ) $
%(inner hom functors) to be a weak equivalence in $\mm$ for all
%covering sieves $ F \rarr C$ of objects $C \in \cc$. This gives a
%model category ``model theoretically enriched'' over $\ttt$.

These examples of model categories of presheaves on a site, as
well as ours, first become interesting after we suitably localize
them. We consider localizations in Section \ref{sec:stalkwise}.
\end{rem}

%\begin{rem}
%The name projective model structure on the category of chain
%complexes of $R$-modules comes from the fact that the cofibrant
%objets are exactly colimits of chain complexes of bounded below
%chain complexes of projective $R$-modules.
%\end{rem}

\section{Tensor structures} \label{sec:tensor}
 The category $\ch ( \nn )$ is a symmetric
 closed tensor category. Let
 $\otimes_R$, or simply $\otimes$, denote the tensor
 product in $\ch (\nn )$. Let $F_R$, or simply $F$, denote
  the inner hom functor in $\ch (\nn )$.
  Our discussion of tensor model categories follows \cite{ss}.
The next Lemma says that all cofibrant objects in the projective
model structure are flat chain complexes.
% (also called $K$-flat
%complexes to distinguish them from complexes of flat modules).
\begin{lem} \label{lem:tensormodel}
Let $K $ be a cofibrant object in $\ch (\nn )$, and let $f \colon X
\rarr Y$ be a weak equivalence. Then $ K \otimes f$ is also a weak
equivalence in $\ch (\nn )$. \end{lem}
\begin{proof} The complex $K$ is a retract of a relative $I$-cell
complex $K'$. Let $ f \colon X \rarr Y $ be a map of presheaves.
There is a natural isomorphism
\[ ( R_C \otimes_R X) (D) = \oplus_{\cc ( D , C ) }
 X (D ) \] for any
two objects $C $ and $ D $ in $\cc$. Hence $K' (D)$ is a directed
colimit of bounded below complex of free $R(D)$-modules.
 The tensor product of a homology-isomorphism
  with a bounded below complexes of free modules
  is again a homology-isomorphism \cite[3.2, 5.8]{spal}.
  Homology commutes with directed
  colimits. Hence
 $ K '\otimes_R X \rarr K' \otimes_R Y$ is a
homology-isomorphism. A retract of a homology-isomorphism is again a
homology-isomorphism, so $ K \otimes f$ is a homology-isomorphism.
 %\[ H_n ( K
%\otimes_R X ) (D) \cong K (D) \otimes_{R(D)} H_n ( X ) ( D) . \]
\end{proof}

The unit object for the tensor product on $\ch (\nn)$ is the chain
complex with $R$ in degree 0. We denote this chain complex by $R$ or
$R[0]$.

\begin{lem} \label{lem:cofibrantunit}
The unit object $R$ is cofibrant. More generally, $R_C$ is
cofibrant, for $C \in \cc$.
\end{lem}

\begin{proof}
The map of chain complexes $0 \rarr R_C$ is the pushout of $i_{U,
-1}$ along the map to the zero chain complex. Hence it is a cofibration.
%acyclic cofibration generator $j_{C ,0} $ in Definition
%\ref{defn:IJ}. Alternatively, one can verify that the map $0 \rarr
%R_C$ has the left lifting property with respect to all acyclic
%fibrations in the projective model structure in Theorem
%\ref{thm:projectivemodelstructure}.
\end{proof}

 Let
$f_1 \colon X_1 \rarr Y_1 $ and $f_2 \colon X_2 \rarr Y_2 $ be two
maps. Then the pushout-product map is the canonical map \[M( f_1
,f_2) \colon colim (Y_1 \otimes X_2 \stackrel{f_1 \otimes
1}{\longleftarrow} X_1 \otimes X_2 \stackrel{1 \otimes f_2
}{\longrightarrow} X_1 \otimes Y_2 ) \longrightarrow Y_1 \otimes Y_2
.\]
 \begin{defn} A model category with a tensor product satisfies the
\mdfn{pushout-product axiom} if $M( f_1 ,f_2)$ is a cofibration
whenever $ f_1 $ and $f_2$ are cofibrations, and $M( f_1 ,f_2)$ is
an acyclic cofibration if $f_1$ or $f_2$ in addition is a weak
equivalence. \end{defn}

\begin{lem} \label{lem:pushoutproduct} The projective model
structure on $\ch (\nn )$ satisfies the pushout-product axiom.
\end{lem}
\begin{proof} The (acyclic) cofibrations are closed under retracts,
transfinite compositions, and pushout. So it suffices to show that
if $f_1$ and $f_2$ are maps in $I$, then $ M( f_1 ,f_2)$ is a
relative $I$-cell complex, and if $f_1 $ is a map in $I$ and $f_2$
is a map in $J$, then $ M( f_1 ,f_2)$ is a relative $J$-cell
complex. Note that $R_{C_1} \otimes R_{C_2}$ is isomorphic to $
R_{C_1 \times C_2 }$, for objects $ C_1 , C_2 \in \cc$. Denote this
object by $R_{12} $ for brevity.
 We have that $ M ( i_{C_1 , 0} ,
i_{C_2 , 0 } ) $ is the inclusion map
 \[ (\cdots \rarr 0  \rarr 0 \rarr R_{12} \oplus R_{12}
 \stackrel{a}{\rarr} R_{12}
 \rarr 0 \rarr \cdots ) \longrightarrow \] \[
 ( \cdots \rarr 0\rarr R_{12} \stackrel{b}{\rarr}
  R_{12} \oplus R_{12} \stackrel{a}{\rarr} R_{12} \rarr 0
  \rarr \cdots
   ) \]
where $ a$ is the fold map and $b$ is given by the identity map on
the first factor and minus the identity map on the second factor.
This is a relative $I$-cell complex. Similarly, the map $ M ( i_{C_1 , 0} , j_{C_2 , 0 } ) $ is a relative $J$-cell complex.
%The map $M ( i_n , i_m )$ is injective with a cofibrant kernel by
%\ref{lem:cibratant}. Hence it is a cofibration by Lemma
%\ref{lem:cofibrations}. The map $M ( i_n , j_0 )$ is a cofibration
%between acyclic objects.
\end{proof}

%\begin{pro} The projective model structure on
%$\ch (\nn )$ is a tensor model category. \end{pro}
%\begin{proof}
% This follows from Lemma
 % \ref{lem:tensormodel} and
% \ref{lem:pushoutproduct}.
%\end{proof}

\begin{defn}
A model category with a tensor product satisfies the \mdfn{monoid
axiom} if $ j \otimes X$ is a weak equivalence for every acyclic
cofibration $j$ and any object $X$, and if pushout and transfinite directed
composition of such maps is again a weak equivalence \cite{ss}. \end{defn}
\begin{lem} \label{monoidaxiom}
 The category $\ch (\nn )$ with the projective tensor model structure
  satisfies
the monoid axiom. \end{lem}
\begin{proof} It suffices to
prove this when $j$ is a map in the set $J$ of generators for the
acyclic cofibrations. This follows since there is a natural
isomorphism
\[ ( R_C \otimes_R X) (D) = \oplus_{\cc ( D , C ) }
 X (D ) \] for any
two objects $C $ and $ D $ in $\cc$. \end{proof}

\begin{defn} \label{defn:a-alg}
 Let $A$ be a monoid in $\ch (\nn )$, i.e.~$A$ is a
differential graded $R$-algebra. Denote  the
category of $A$-modules in $\ch (\nn )$ by \mdfn{$A $-$\ch (\nn )$}. \end{defn}
%There is a
%forgetful functor $k \colon A $-$ \ch (\nn) \rarr \ch (\nn) $. It
%has a left adjoint functor given by tensoring from the left by $A$
%over $R$. The adjunction $ ( A \otimes_R - , k)$ gives that the
%sources of the maps in $\Sigma_A I $ and $ \Sigma_A J$ are small,
%and Lemma \ref{monoidaxiom} gives that $\Sigma_A J$-cell complexes
%are in $W$. Hence
\begin{lem} The category $A $-$\ch (\nn )$ inherits a cofibrantly
generated model structure from $\ch (\nn )$ via the forgetful
functor $A $-$ \ch (\nn) \rarr \ch (\nn) $. If $A$ is a symmetric
monoid, then $A $-$\ch (\nn )$ is a tensor category
(with tensor product over $A$)
 and the model structure on $A $-$\ch (\nn )$ satisfies the
pushout-product axiom and the monoid axiom. \end{lem}
\begin{proof} This follows from
Lemmas \ref{lem:tensormodel}, \ref{lem:pushoutproduct}, and
 \ref{monoidaxiom} and \cite[3.1]{ss}.
\end{proof}
\begin{lem}
Let $A$ be a symmetric monoid. Then the category of $A$-algebras
inherits a model structure via the forgetful functor from
$A$-algebras to $A $-$ \ch ( \nn ) $.
\end{lem}
\begin{proof} This follows from
 Lemmas \ref{lem:tensormodel}, \ref{lem:cofibrantunit},
  \ref{lem:pushoutproduct}, and
 \ref{monoidaxiom} and \cite[3.1]{ss}.
\end{proof}
%We also have that the
%category of symmetric $A$-algebras inherits a model structure from
%the forgetful functor from symmetric $A$-algebras to $A $-$ \ch (\nn
%)$.

\section{Quasi-simplicial model structures} \label{sec:quasi}
We model theoretically enrich $ \ch (\nn )$ in simplicial sets.
%One might argue that it would
%be more natural to model theoretically enrich $ \ch (\nn)$, with the
%projective model structure, in the model category of unbounded chain
%complexes of abelian groups with the projective model structure. (Or
%over the category itself.) We choose not to do so.
%I am grateful to Brooke Shipley for explaining to me why the obvious
%guess of a simplicial structure does not quite work \cite[2.7]{ss2}.
  %This section is not essential for the paper.
  %We first enrich $\ch (\nn)$ over
% itself. We give $ \ch (\nn)$ a tensor and cotensor structure over
 % itself from
  %the tensor product and inner hom functor. The inner hom functor is
  %also the mapping complex. This is a $\ch (\nn )$-model category
  %in the sense of \cite[9.1.2]{hir}
  % (with simplicial sets replaced by $\ch (\nn)$). This follows
  % from Lemma \ref{lem:pushoutproduct} and the proof of Theorem
 % \ref{thm:stalkwise}.
%This is true for both the projective and the stalkwise model
%structures (when we use the same model structure on $\ch ( \nn)$ as
%on its enrichment).
We define  a weakening of the axioms for a simplicial model category
\cite[9.1.6]{hir}.

\begin{defn} \label{def:quasisimplicial} A
\mdfn{quasi-simplicial} model category is a model category $\kk$
that is a simplicial category (see \cite[9.1.2]{hir})
  satisfying the axioms below.
   Map denotes the (based) simplicial mapping
space, and $X \square S $ and $F_{\square} ( S , X )$ are the tensor
and cotensor of $X \in \kk$ by a simplicial set $S$, respectively.
\begin{description}
 \item[weakM6]

Let $\sss$ be the category of simplicial sets.
Let $X , Y $ be objects in $\kk$ and let $S $ be an object in $\sss$.
 There is a natural isomorphism of simplicial
sets
\[ \text{Map} ( X \square S , Y ) \cong \text{Map}
(X , F_{\square} ( S , X )). \] There is a natural isomorphism of \emph{sets} \[ \sss
( S , \text{Map} ( X , Y ) ) \cong \kk ( X \square S , Y ) .\]
\item[M7] Let $i\colon A \rarr B$ be a cofibration in $\kk$
and $ f \colon
 X \rarr Y$ a
fibration in $\kk$. Then the map \[ i^* \times f_* \colon \text{Map}
(B ,X ) \rarr \text{Map} (A ,X ) \times_{\text{Map} (A ,Y )}
\text{Map} (B ,Y ) \] is a fibration of simplicial sets. If, in
addition, $i$ or $f$ is a weak equivalence, then $i^* \times f_*$ is
a weak equivalence.
%\item Let $f : X \rarr Y $ be a cofibration in $\kk$
%and $j : S \rarr L$ is an injection of simplicial sets. Then the
%pushout-product map
%\[M( f ,j) \colon
%colim (Y \square S \stackrel{f \otimes 1}{\longleftarrow} X \square
%S \stackrel{1 \otimes j}{\longrightarrow} X \square L)
%\longrightarrow Y \square L\] is a cofibration in $\kk$.
 %Moreover, if in addition
%$ f $ or $j$ is a weak equivalence, then $M( f, j )$ is an acyclic
%cofibration.

% The natural map
% $ (X \square K )\square L \rarr X \square ( K \times L)$
%is a weak equivalence, whenever $X$ is a cofibrant object in $\kk$
%and $K $ and $L$ are arbitrary simplicial sets. (The map corresponds
%to the identity map $ 1_{X \square (K \times L)}$ under a sequence
%of adjunctions using the isomorphism $ \text{Map} ( K , \text{Map} (
%L , X \square ( K \times L )) \cong \text{Map} ( K \times L , X
%\square ( K \times L ))$ that is part of the simplicial structure.)
\end{description}
\end{defn}

\begin{lem} \label{M8} Let $\kk$ be a quasi-simplicial model category. Then the  natural map
 $ X \square * \to X$ (corresponding to $1_X$ under the adjunction in weakM6)
 is an isomorphism, for all $X \in \kk$. \end{lem}
 \begin{proof}
 The Yoneda lemma applied to the composition of isomorphisms
 \[ \kk ( X , Y ) \cong \text{Map} ( X , Y)_0 \cong \sss ( \ast , \text{Map} (X , Y) ) \cong \kk ( X \square * , Y )  \] gives the result.
  \end{proof}

The second axiom has an
 equivalent formulations in terms of the tensor or cotensor
 functors instead of the simplicial mapping space
  \cite[9.3.6]{hir}. One implication of M7 is that $ X \square S
  \rarr X \square T$ is a weak equivalence in $\kk$ whenever $ X$ is
  cofibrant in $\kk$ and $S \rarr T$ is a weak equivalence of
  simplicial sets. Combined with Lemma \ref{M8} this gives that the suspension
of a cofibrant object $X$ is equivalent to $X \square S^1$ and the
loop of a fibrant object $Y$ is equivalent to $ F_{\square} ( S^1 ,
Y)$.

In a simplicial structure the last isomorphism in weakM6 is an isomorphism of simplicial sets.
In a quasi-simplicial structure, unlike a simplicial structure, repeated applications of the tensor and cotensor
functors need not respect the
 cartesian product (based: smash product) on simplicial sets \cite[9.1.11]{hir}.

\begin{lem} \label{simplicial}
The projective model structures on $\ch (\nn )$ is quasi-simplicial.
\end{lem}

\begin{proof}
We make use of the Dold-Kan adjunction \cite[8.4]{wei}. We note that
the inclusion functor, \mdfn{$k$}, from nonnegative chain complexes to
unbounded chain complexes is left adjoint to the truncation functor,
\mdfn{$\tau_{\geq 0}$}, given by sending \[\cdots \rarr X_2 \rarr X_1 \rarr
X_0 \rarr X_{-1} \rarr X_{-2} \rarr \cdots \] to \[\cdots \rarr X_2
\rarr X_1 \rarr \ker (X_0 \rarr X_{-1} ) \rarr 0 \rarr \cdots . \]

 Let \mdfn{$D$}
denote the right adjoint of the normalized chain complex functor
\mdfn{$N$}. Let $S$ be a simplicial set. Let $R [ S]$ be the free
simplicial $R$-module, and let $N R [S]$ be the associated chain
complex. We define the simplicial tensor to be $ X \otimes_R k N R
[S]$ and the simplicial cotensor to be $ F ( k N R [S] , X )$, for
$X \in \ch ( \nn )$. The simplicial hom functor, $\text{map} ( X , Y
)$, is defined to be
\[ D j \tau_{\geq 0} \Gamma (F ( X, Y)) ,\] for $X , Y \in \ch (\nn )$, where
\mdfn{$\Gamma$} denotes the global sections functor, \mdfn{$j$} is the forgetful
functor from the category of $\Gamma R $-chain complexes to the
category of chain complexes of abelian groups. Clearly the
adjunctions in
%Definition \ref{def:quasisimplicial}
  axiom weakM6 are satisfied.

 Note first that the weakened version of
 axiom M6 still implies that axiom M7 has an alternative
 formulation as a pushout-product axiom in terms of the tensor
 \cite[9.3.7]{hir}.
  If $i \colon S \rarr S'$ is an
inclusion of simplicial sets, then $kNR [ i] : kNR [S] \rarr kNR
[S']$ is a cofibration in the projective model structure, and if $i$
is a weak equivalence, then $kNR [ i] $ is a presheaf homology
isomorphism in $\ch ( \nn )$. Hence axiom M7 follows from
    Lemma \ref{lem:pushoutproduct}.
    %Axiom M8 is satisfied \cite[2.2,4.2]{ss2}.
% since $k N R [j]$ is a
% cofibration in $\ch (\nn) $ when $j$ is an inclusion of simplicial
% sets, and an acyclic cofibration when $j$ is an injective weak
% equivalence of simplicial sets.
\end{proof}

\section{The stalkwise model structure }
\label{sec:stalkwise}

Given a (skeletally) small Grothendieck site $\cc$. The following
model structure is called the \mdfn{stalkwise model structure} on
$\ch (\nn)$.

\begin{thm} \label{thm:stalkwise}
There is a proper quasi-simplicial
  cofibrantly generated   model structure
 on $\ch (\nn )$. The  weak
equivalences are stalkwise homology-isomorphisms and the cofibrations are retracts of relative $I$-cell complexes. With this model structure the tensor closed structure  on $\ch ( \nn )$
satisfies the pushout-product axiom and the monoid axiom.
\end{thm}
\begin{proof}
The stalk functors from $\ch (\nn )$ to the categories of chain
complexes of $R_p$-modules respect both pushout and pullback
squares, levelwise surjective maps, and levelwise injective maps.
Furthermore, it takes homology-isomorphism of presheaves to
homology-isomorphisms of $R_p$-modules.

The functor $ (H_n )_p $, for each point $p$ in $\ee$,
  is a homology theory in $ \ch (\nn )$
 which commutes with arbitrary directed colimits.
Hence
 we can Bousfield localize $\ch ( \nn )$ with respect to $ (H_n )_p
$-equivalences for $n \in \zz$ and $p$ in $\pt$. This requires a
cardinality argument involving the cardinality of $ R (C)$ for all
$C \in \cc$. More  details are included in Appendix
\ref{appendix:bousfield}.

The model structure is proper and satisfies the pushout-product
axiom and the monoid axiom. This follows from
  Theorem \ref{thm:projectivemodelstructure}, Lemmas
\ref{lem:pushoutproduct} and \ref{monoidaxiom}, and stalkwise
verification of weak equivalences. The model structure is
quasi-simplicial, with a quasi-simplicial structure given by the
quasi-simplicial structure on the projective model structure (Lemma \ref{simplicial}). This
follows because the stalkwise model structure satisfies the pushout-product axiom.
\end{proof}

%[[The author does not have a simple description of the fibrations
%and the fibrant objects in the stalkwise model structure on $\ch
%(\nn )$.
%One might start with a bounded below complex and try to build a
%fibrant replacement inductively. It would be interesting to find a
%nice set of acyclic cofibrant generators. We have that the
%``skyscraper presheaf'' functor, with respect to any point in $\ee$,
%applied to any levelwise surjective chain complex map is a fibration
%in projective model structure on $\ch ( \nn)$. ]]

\begin{pro}
Assume that $L^2 R_C$ is small in $\mm$, for all $C \in \cc$.
Then there is a proper quasi-simplicial cofibrantly generated
model structure on $\ch (\mm )$, with cofibrant generators $L^2 I$
and acyclic cofibrant generators $L^2 J$. The weak equivalences are
the stalkwise homology-isomorphisms and the cofibrations are
retracts of relative $L^2 I$-cell complexes. The model structure
satisfies the pushout-product axiom and the monoid axiom.
\end{pro}

\begin{proof}
By our assumptions the sources of $ L^2 I$ and $L^2 J$ are small. Relative $L^2 J$-cell
complexes in $\ch (\mm)$ are stalkwise homology-isomorphisms. Hence
the result follows by Theorem \ref{thm:stalkwise} and
\cite[11.3.2]{hir}.
The model structure is proper, quasi-simplicial, and satisfies the pushout-product axiom and the monoid axiom. This is verified as in the proof of Theorem \ref{thm:stalkwise}
%Since $L^2$ respects both colimits and pullback
%squares the resulting model category is proper and satisfies both
%the pushout-product axiom and the monoid axiom \cite[II.4.1]{sga4.1}.
\end{proof}
If $X$ is a Noetherian topological space and $R$ is a ring of sheaves on $X$, then $ L^2 R_U$ is small for all open subsets $U$ of $X$ (since direct sums of sheaves in the category of presheaves are themselves sheaves).
%\begin{proof}
%Directed colimts commute with $i\colon \mm \rarr \nn$? That $L^2
%R_C
%$ is small?
%The class of stalkwise equivalences is closed under pushout along a
%cofibration and under pullback along fibrations. This follows since
%we can check this on stalks.
%\end{proof}

\begin{pro} \label{pro:sheaf-presheafQEqual}
Assume that $L^2 R_C$ is small in $\mm$, for all $C \in \cc$.
Let $\ch ( \mm )$ and $ \ch (\nn )$ both have the stalkwise model
structure. The map of topos $ i \colon \mm \rarr \nn$ and $L^2
\colon \nn \rarr \mm$ induces a Quillen equivalence of model
categories $ i \colon \ch ( \mm) \rarr \ch (\nn)$.
\end{pro}
\begin{proof} The sheafification functor respects cofibrations and weak equivalences.
%Let $X$ be a preasheaf and let $Y$ be a sheaf.
A map $   L^2 X \rarr Y$ is a
stalkwise homology-isomorphism if and only if the adjoint map $ X
\rarr i (Y)$ is a stalkwise homology-isomorphism for every $X \in \ch
( \nn )$ and $Y \in \ch ( \mm )$.
\end{proof}

%\begin{rem} The points of the topos of $\cc$-sheaves for the chaotic
%topology (only the identity maps are covers) are the evaluation
%functors at $C$, for $C \in \cc$ \cite[IV.6.8.6]{sga4.1}. So the %projective and stalkwise
%model structures on $\ch (\nn )$ coincide. Hence the projective model
%structure is a special case of the stalkwise model structure.
%\end{rem}
Recall Definition \ref{defn:derivedcat} of the derived category
$\dd_R$.
%, of chain complexes of sheaves of $R$-modules on the site $\cc$.
\begin{pro} Assume that the topos $\shc$ has enough
points. The homotopy category of $\ch ( \nn )$ with the stalkwise model structure is  equivalent to $\dd_R$ as tensor
triangulated categories.
\end{pro}
\begin{proof} Since $\ee$ has enough points the
  stalkwise homology-isomorphism and
the sheaf homology-isomorphism coincide \cite[IV.6.4.1]{sga4.1}.
%The result follows from Proposition \ref{pro:sheaf-presheafQEqual}.
\end{proof}

%We consider some full subcategories of $\ch (\nn )$. These full
%subcategories need not be complete nor cocomplete.

%\begin{pro} The full subcategory of $\ch (\nn )$ of objects whose
%homology (stalkwise homology) is bounded below is a full model
%category of $\ch (\nn )$ with the projective model structure (the
%stalkwise model structure). \end{pro}

Let $\ch ( \nn)_{+} $ be the full subcategory of $\ch (\nn)$
  consisting of chain complexes $\{ X_n \} $ such that $X_n = 0$, whenever $n <0$.
The acyclic cofibrant generators $J'$ is the set of maps $j_{C,n} $,
for $n \geq 0$ and $C \in \cc$, and the cofibrant generators $I'$ is the set of maps
$ i_{C,n}$, for $n \geq 0$, together with the maps $ 0 \rarr R_C
[0]$, for $C \in \cc$.
\begin{pro} \label{pro:chainsgeq0}
There is a cofibrantly generated model structure on $\ch (
\nn)_{+} $ such that the weak equivalences are presheaf
homology isomorphisms. The cofibrations are retracts of relative
$I'$-cell complexes. The fibrations are maps $f \colon X \rarr Y$
such that $ f_n (C) $ is surjective, for all $C \in \cc$ and $n \geq
1$. \end{pro}

\begin{proof} This is essentially contained in the proof of Theorem
\ref{thm:projectivemodelstructure}. Note that the class
$\text{inj}\, ( 0 \rarr R_C [0]) $ in $\ch ( \nn)_{+} $
consists of maps that are surjective in degree 0 when evaluated at
$C$. \end{proof}

\begin{pro}
There is a cofibrantly generated model structure on $\ch (
\nn)_{+} $ such that the weak equivalences are stalkwise
homology isomorphisms, and the cofibrations are retracts of relative
$I'$-cell complexes.
%The cofibrant generators are $ i_{C,n}$, for
%$n \geq 0$ together with the maps $ 0 \rarr R_C [0]$, for $C \in
%\cc$.
\end{pro}

\begin{proof} This follows by localizing the model structure in
 Proposition
\ref{pro:chainsgeq0} with respect to stalkwise equivalences. Note
that $\ch ( \nn)_{+} $ is closed under directed colimits.
\end{proof}

\section{Some t-structures on derived categories}
\label{sec:t-structures}

 We construct t-structures on the homotopy
category of $\ch (\nn )$ with the stalkwise model
structures. These t-structures interact
well with the model structure on $\ch (\nn)$. More precisely, they
all arise from a t-model structure on $\ch (\nn )$.
% For the definition of
% a t-structure on a triangulated category see \cite[2.1]{tfi}.
Homological grading of t-structures is used. So $ \dd_{\geq n} ,
\dd_{\leq n -1 } $ corresponds to $ \dd^{\leq -n} , \dd^{\geq -n +1
} $ in cohomological notation.
Let $\kk$ be a stable model category, with a t-structure on its triangulated homotopy category $\text{Ho} (\kk) $ \cite[7.1]{hov}.
\begin{defn}
The class of
\mdfn{(co-)$n$-equivalences} in $\kk$ is the class of maps $f$ in $\kk$ such
that the homotopy type of $\hofib (f)$ is in $\text{Ho}(\kk)_{\geq n}$
($\text{Ho}(\kk)_{\leq n-1}$) \cite[3.1]{tfi}.
\end{defn}
We now make precise what we mean by lifting a t-structure on $\text{Ho} (\kk)$ to $\kk$.

\begin{defn} A (weak) \mdfn{t-model category}
 is a proper quasi-simplicial
stable model category $\kk$ with functorial factorizations equipped
with a t-structure on its homotopy category together with a
functorial factorization of maps in $\kk$ into $n$-equivalences
followed by co-$n$-equivalences, for each $n \in \zz$.
\end{defn}
This is a weakening of the definition of a t-model structure in
\cite[4.1]{tfi}. We require the model structure to be
 quasi-simplicial instead of simplicial. This is a harmless
 weakening and all the results of \cite{tfi} are still valid
 (with simplicial replaced by quasi-simplicial).

We are mainly interested in t-structures on $\ch (\nn )$ with the
stalkwise model structure, but we consider a more general framework.
Let $\g$ be a proper quasi-simplicial cofibrantly generated stable cellular
model category together with a t-structure on its homotopy category.
Let $I$ be a set of cofibrant generators of $\g$ \cite[12]{hir}. We make the following assumptions:
\begin{enumerate} \item the maps in $I$ have small sources; \item
the heart of the t-structure is the category of sheaves, $ \mm $,
(or of presheaves, $\nn$,) of $R$-modules for a ring of sheaves $R$
on a Grothendieck site $\cc$; and \item the heart functor, $ \hh
\colon \g \rarr \mm$ is $\sigma$-uniform, for a cardinal number $\sigma$ (see Definition
\ref{defn:sigmauniform}), respects sums, and directed colimits of
relative $I$-cell complexes.
\end{enumerate} The topos  \mdfn{$\ee$} is assumed to have a set, $\pt$, of
isomorphism classes of points. Let $\dd$ denote the homotopy category of $\g$.
 Let
\[ \mdfn{d} \colon \pt \rarr \zz \cup \{ \pm \infty \} \]
be a function.
 We construct t-model structures on $\kk$ by shifting the original t-structure such that at each point of $\ee$ in the isomorphism class  $p \in \pt$ the shift is given  by $\Sigma^{d (p)}$.
%     isomorphism class of points $p \in \pt $.

\begin{pro} \label{t-structures}
 There is a t-model structure (with simplicial relaxed to
quasi-simplicial)
   on $\g$ such that \[ \dd_{\geq 0 } = \{ X \, | \, (
\hh_{n_p} (X ) )_p = 0 \text{ for all } n_p < d (p) \} .\]
\end{pro}
\begin{proof} The associated class
of
  $n$-equivalences, $W_n$, consists of all maps $f $ such that
$( \hh_{n_p } (f) )_p$ is an isomorphism for all $n_p < d (p)+n $
and $( \hh_{d (p)+n} (f) )_p$ is a surjection if $|d (p)| < \infty$.

% Let denote the
%class of $n$-equivalences.
 The pushout of a
$W_n$-map along a cofibration is again a $W_n$-map. For each point
$p$ in the topos $\shc$ the functor
 $ (\hh )_p$ from
$\g$ to abelian groups respects sums and directed colimits of
relative $I$-cell complexes. We can localize the category $\g$ with
respect to the $\zz$-indexed homology theory whose $n$-th functor is
\[ X \mapsto \textstyle\bigoplus_p ( \hh_{d (p) + n } (X) )_p \]
using Proposition \ref{pro:bousfield}. See also
  \cite[7.5]{tfi}. \end{proof}

The full subcategory $\dd_{\leq -1}$ is given by \[ \{ Y \in \dd \ |
\ \dd ( X , Y ) =0 \text{ for } X \in \dd_{\geq 0} \}. \] We
describe $\dd_{\leq -1}$ more explicitly in Section \ref{sec:ddgeq1}
when $\g $ is $\ch (\nn)$.

%Two different functions $d$ might give rise to the same t-structure
%on $\g$.

\begin{cor} \label{cor:localization}
Let $Z$ be a subset of $ \pt $. Then there is a proper
quasi-simplicial model structure on $\g$ such that the weak
equivalences are $\oplus_{p \not\in Z} (\hh_* )_p $-isomorphisms and
the cofibrations are retracts of relative $I$-cell complexes.
\end{cor}
\begin{proof}
Let $d_Z \colon \pt \rarr \zz \cup \{ \pm \infty \} $ be the
function defined by letting $ d_Z ( z ) = \infty $, for $ z \not\in
Z$, and $ d_Z ( z ) = -\infty $, for $ z \in Z$. The result follows
from Proposition \ref{t-structures} applied to the function $d_Z$.
\end{proof}

We refine the t-structure in Proposition \ref{t-structures} by
taking the structure of the ring $R$ into account. Let $ d \colon
\coprod_{p \in \pt } \spec R_p \rarr \zz \cup \{ \pm \infty \} $ be
a function. We can localize the category $\g$ with respect to the
$\zz$-indexed homology theory whose $n$-th functor is
\[ X \mapsto \oplus_{p,\mathfrak{p}}
(( \hh_{d ({p , \mathfrak{p} })+ n } (X) )_p)_{\mathfrak{p}}
  .\]

\begin{pro} \label{t-structuresrefined}
There is a t-model structure (with simplicial relaxed to
quasi-simplicial)
%(with simplicial relaxed to quasi-simplicial)
 on $\g$ such that
\[ \dd_{\geq 0 } = \{ X \, | \, (
 \hh_{n_{p,\mathfrak{p}}} (f) )_p)_{\mathfrak{p}}
 = 0 \text{ for all } n_{p,\mathfrak{p}} < d ({p , \mathfrak{p} })
 \} .\]
\end{pro}
\begin{proof}
The corresponding class of
 $n$-equivalences
 %, $W_n (d) $,
 is the
  class of maps $f$ in $\g$ such that
 \[ ((\hh_{n_{p,\mathfrak{p}}} (f) )_p)_{\mathfrak{p}}\] is an
isomorphism for all $ n_{p,\mathfrak{p}} < d ({p , \mathfrak{p} })+
n $, and
%$((\hh_{d_{p,\mathfrak{p}} +n } (f) )_p)_{\mathfrak{p}}$ is
 a
surjection for $ n_{p,\mathfrak{p}} = d ({p , \mathfrak{p} }) + n $
if $|d ({p , \mathfrak{p} })| < \infty$. The result follows from
Proposition \ref{pro:bousfield} since the class of $n$-equivalences
is closed under pushouts along cofibrations and the sources of the
cofibrant generators are small.
\end{proof}

If $d$ is constant on each $\spec R_p$, then Proposition
\ref{t-structuresrefined} reduces to Proposition \ref{t-structures}.
%We extend the class of t-structures on the homotopy category $\dd$
%of $\g $. We are not able to show that this set of t-structures
%comes from a t-model structure on $\g $.
We consider another t-model structure on $\g$. Let \[ d \colon
\coprod_{p \in \pt } \spec R_p \rarr \zz \cup \{ \pm \infty \} \] be
a function.
\begin{pro} \label{general-t-structures}
 There is a t-model structure (with simplicial relaxed to
quasi-simplicial) on $\g$ such that \[ \dd_{\geq 0 } = \{ X \, | \,
((
 \hh_{n_{d , p ,\mathfrak{p}} } (X ) )_p)_{\mathfrak{p}}
\otimes_{R_{\mathfrak{p}}} R_{\mathfrak{p}} / \mathfrak{p} = 0
\text{ for all } n_{d , p ,\mathfrak{p}} < d ( p ,\mathfrak{p} ) \}
.\]
\end{pro}
\begin{proof} Note that $ - \otimes_{R_{\mathfrak{p}}}
R_{\mathfrak{p}} / \mathfrak{p} $ is exact on the category of $
R_{\mathfrak{p}}$-modules since $ R_{\mathfrak{p}} / \mathfrak{p} $
is a field. The corresponding class of
 $n$-equivalences is the class of maps  $f \colon X \rarr Y $ such that
\[ (( \hh_{n_{d , p ,\mathfrak{p}} } (f ) )_p)_{\mathfrak{p}}
\otimes_{R_{\mathfrak{p}}} R_{\mathfrak{p}} / \mathfrak{p} \] is an
isomorphism, for $ n_{d , p ,\mathfrak{p}} < d ( p ,\mathfrak{p} )$,
and an surjection for $ n_{d , p ,\mathfrak{p}} = d ( p
,\mathfrak{p} )$ if $ |d ( p ,\mathfrak{p} ) | < \infty$. The
t-model structure is obtained from
%given by localizing with respect to the
%functors:
%\[ X \mapsto \oplus_{p,\mathfrak{p}}
%(( \hh_{d_{p,\mathfrak{p}} +n } (X) )_p)_{\mathfrak{p}}
%\otimes_{R_{\mathfrak{p}}} R_{\mathfrak{p}} / \mathfrak{p}
% . \]
Proposition \ref{pro:bousfield} since the sources of the cofibrant
generators are small and the class of $n$-equivalences is closed
under pushouts along cofibrations.
%Call this functor $K_d$.
  \end{proof}

\begin{exmp}
Propositions \ref{t-structuresrefined} and
\ref{general-t-structures} apply to the category $\ch (\nn )$ with
the stalkwise model structure and the standard t-structure on its
homotopy category. The heart valued homology functor associated to
the standard t-structure on the derived category $\dd$ is the usual
homology of a chain complex. It respects sums and directed colimits
of relative $I$-cell complexes.

Let $A$ be a monoid in $ \ch (\nn )$ such that the sheaf valued
homology is zero in negative degrees. The assumption on $A$ gives
that there is  a standard t-structure on the homotopy category of $A$-$ \ch (\nn )$.
 Propositions \ref{t-structuresrefined} and
\ref{general-t-structures} apply to the category of $A $-modules in
$ \ch (\nn )$ with the standard t-structure. See Definition
\ref{defn:a-alg}.
 \end{exmp}

A class of t-structures on $\ch (\nn )$ with the
projective model structure can be constructed using  another technique. Let $I'$ be
a subset of $I$ such that whenever $ i_{C,n} \in I'$, then $ i_{C
,m} \in I'$ for all $m \leq n$. We associate to $I'$ a class of maps
closed under (nonnegative) suspensions: \[ W (I') = \{ f \, | \,
H_{m-1} (f) (C) \text{ is an isomorphisms in } \nn \text{ and } \]
\[ \, H_{m} (f) (C) \text{ is a surjection in }
\nn, \text{ whenever } i_{C,m} \in I' \} . \]

\begin{lem} There is a t-model
structure on $\ch (\nn )$ with the projective model structure such
that the associated class of $0$-equivalences is $W ( I')$.
\end{lem}
\begin{proof} This follows from the proof of Theorem
\ref{thm:projectivemodelstructure} and \cite[4.13]{tfi}.
\end{proof}

\section{A description of $\dd_{\leq 0}$ for the derived category}
\label{sec:ddgeq1}
 Recall that if $\ee$ has enough points, then the homotopy category
 of
 $ \ch ( \nn)$ with the
stalkwise model structure is equivalent to the derived category,
$\dd $.
 We describe the full subcategory,
$\dd_{\leq 0}$, of $\dd $ for t-structures obtained from functions $d
\colon \pt \rarr \zz \cup \{ \pm \infty \}$ such that $ d^{-1} ( [
n, \infty ] ) $ is an open subset of $\pt $, for every $ n \in \zz
\cup \{\pm \infty \}$. We first recall some terminology and a Lemma.

The support of an object $ X$ in $\ee$ is defined to be \[ \sup (X)
= \{ p \in \pt \ | \ X_p \not= \emptyset \} . \] Note that if $X
\rarr Y$ is a map in $\ee$, then $\sup (X) \subset \sup (Y)$. Recall
that the topology on $\pt $ is generated by a basis of open sets
consisting of $ \sup ( S)$, for all subobjects $S$ of the terminal
object $\bullet $ of $\ee$ \cite[IV.7.1.7]{sga4.1}. Neighborhoods of
points are defined in \cite[IV.6.8]{sga4.1}.

\begin{lem} \label{lem:supportopen}
Let $X $ be an object in $\ee$ and let $p $ be a point in $\ee$.
Suppose given an element $x \in X_p$ and an open subset $U$ of $\pt$
containing $p$. Then there exists
  an object $C \in
 \cc$ with support contained in $U$ and a map $ C \rarr X$
 in $\ee$ such
 that $x$ is in
 the image of $ C_p \rarr X_p$. \end{lem}
In other words, $p$ has a neighborhood with support in $U$.
\begin{proof} By definition of the topology on $\pt$ there is a
subobject $S$ of $\bullet$ such that $p \in \sup (S) \subset U$. The
stalk, $ S_p = \bullet$, is the colimit of $S(C)$ for neighborhoods
$(C , c\in C_p )$ of $p$. Hence there exists a neighborhood $(C ,
c\in C_p )$ of $p$ such that  $x$ in the image of $ C_p \rarr X_p$ and $S (C) \not= \emptyset$ (so there is a
map $ C \rarr S$ in $\ee$). Since $\sup (C) \subset \sup (S)$ the
claim follows.
 \end{proof}

\begin{pro} \label{pro:tBBD} Assume that $\ee$ has enough points.
Let $d \colon \pt \rarr \zz \cup \{\pm \infty \}$ be a function such
that $ d^{-1} ( [ n, \infty ] ) $ is an open subset of $\pt $, for
every $ n \in \zz \cup \{\pm \infty \}$. Then there is a t-model
structure on $\ch (\nn )$ (with simplicial relaxed by
quasi-simplicial) such that
\[ \dd_{\geq 0 } = \{ X \, | \, (
\hh_{n_p} (X ) )_p = 0 \text{ for all } n_p < d (p) \} \] and
\[ \dd_{\leq 0 } = \{ X \, | \, (
\hh_{n_p} (X ) )_p = 0 \text{ for all } n_p > d (p) \} .\]
\end{pro}

\begin{proof} Except for the description of $\dd_{\leq 0}$
 the result follows from Proposition
\ref{t-structures}.

We construct an explicit truncation functor. For each object $C \in
\cc$ let $n (C) \in \zz \cup \{ \pm \infty \}$ be the largest number
such that $n(C) \leq d(p)$, for all $p \in \sup (C)$. Define
$X_{\geq 0} $ to be the complex such that $(X_{\geq 0})_k (C)$ is
$X_k (C)$ for $k> n(C)$, $\ker ( X_{n(C)} (C) \rarr X_{n(C)-1} (C))$
for $k = n(C)$, and $0$ for $k < n(C)$. There is a canonical inclusion
map $
X_{\geq 0} \rarr X$ and we denote the cokernel by $X_{\leq -1}$. The truncation functors give well defined functors in the homotopy category $\dd$.
The
quotient map $ X \rarr X_{\leq -1}$ is a fibration (between fibrant objects) since it is
levelwise surjective. Since  $X_{\geq 0}$ is its fiber
\[X_{\geq 0} \rarr X \rarr X_{\leq -1} \rarr \Sigma X_{\geq 0}
\] gives a triangle in the homotopy category \cite[6.2.6, 6.3, 7.1]{hov}. The map $ X_{\geq
0} \rarr X$ is a stalkwise equivalence if and only if $ X \in \ch
(\nn )_{\geq 0}$, and $ X \rarr X_{\leq -1} $ is a stalkwise
homology isomorphism if and only if $ X \in \dd_{\leq -1}$ by Lemma
\ref{lem:supportopen}. Both $ X \mapsto X_{\geq 0}$ and $ X \mapsto
X_{\leq -1}$ are idempotent functors, and $ (X_{\geq 0} )_{\leq -1}
= (X_{\leq -1} )_{\geq 0} =0$. Hence $ \dd ( X_{\geq 0 } , Y_{\leq
1}) =0$, for all $X , Y \in \dd$.

\end{proof}

\begin{cor} The heart of the t-structure
in Proposition \ref{pro:tBBD} is given by
\[ \{ X \, | \, (
\hh_{n_p} (X ) )_p = 0 \text{ for all } n_p \not= d (p) \} .\]
\end{cor}

\begin{rem} The assumption on $d $ in Proposition \ref{pro:tBBD}
 is not optimal. There might be more general
functions $d$ that have the same description of $\dd_{\leq 0}$. For
example if $ d^{-1} ( [ n , \infty ] )$ is closed (instead of open),
for all $n \in \zz \cup \{ \pm \infty \}$, then Proposition
\ref{pro:tBBD} is still valid (make an explicit construction of the
truncation $X \rarr X_{\leq -1 }$).
%We can also give explicit
%description of $\dd_{\leq -1}$ in Proposition
%\ref{t-structuresrefined} under reasonable assumptions.
\end{rem}

\section{Examples and comparisons} \label{sec:examples}
\begin{exmp}[Rings]
Let $\cc$ be the one morphism site. Then a ringed topos is the
category of sets together with a ring $R$, and $\mm = \nn$ is the
category of $R$-modules. The projective and the stalkwise model
structures on $\ch (\mm)$ coincide. The localization in Corollary
\ref{cor:localization} has been constructed by Neeman for the one
morphism site \cite[3.3]{nee}.

Let $R$ be a Noetherian ring. Stanley has constructed the
t-structures in
  Proposition \ref{t-structuresrefined} on the full
subcategory of $\dd_R$ consisting of complexes whose homology groups
are finitely generated in each degree and bounded above and below
\cite{Sta}. He also shows that there are no other t-structures on
this full subcategory of $\dd$ \cite[5.3]{Sta}.
%This implies that if a t-structure
%constructed in Proposition \ref{general-t-structures} restricts to
%give a t-structure on
% Stanley's full subcategory of $\dd_R$, when $R$ is a Noetherian
% ring, then the restricted t-structure is of the form
% constructed in Proposition
% \ref{t-structuresrefined}.
%Assume that $R_p$ is a Noetherian ring for all $p \in \pt $. Then
%the t-model structures in Propositions \ref{t-structuresrefined} and
%\ref{general-t-structures} restrict to the full subcategory of $\ch
%(\nn )$ consisting of objects whose stalkwise homology is bounded,
%above and below, and finitely generated for all points in $\ee$. In
%particular, if $R$ is a Noetherian ring the t-structures on
%Stanley's full subcategory of $\dd_R$ lift to t-model structures on
%a full model subcategory of $\ch ( \nn )$.
\end{exmp}

\begin{exmp}[Perverse t-structures] We consider the
 category of $R$-modules for a
ringed space $(S, R)$. A space is said to be sober if all closed
irreducible sets have a generic point. Let $k \colon S \rarr
S_{\text{sob}}$ be the universal inclusion into a sober space. The
points of $S_{\text{sob}}$ corresponds to closed irreducible subsets
of $S$, and the map $k$ is given by sending a point $s \in S$ to the
closure of $s$ in $S$. The set of (isomorphism classes of) points of
$(S , R)$ is the space $ S_{\text{sob}}$ \cite[IV.7.1.6]{sga4.1}.

%The set of open subsets of $S$ and $S_{\text{sob}}$ are the same.

%If $S$ is a Noetherian space, then $ R_U$ is a sheaf for every open
%subsets of $S$ (this is the set of isomorphism classes of points in
%this case). We also have that arbitrary sum of sheaves in the
%category of presheaves agree with the sum in the category of
%sheaves.
% Hence we get that the cofibrant complexes in $\ch (\nn )$
%are actually sheaves.

%One can also consider ringed spaces with an action by a discrete
%group $G$.
  We now assume that $S = S_{\text{sob}} $ and compare our t-structures to the
  perverse t-structures introduced by
Beilinson, Bernstein, Deligne, and Gabber \cite{bbd}. They   consider
a nonempty finite  partition $\{ S_a \}_{a \in A}$ of $S $
into locally closed sets, together with  a function $ p \colon A \rarr \zz$,
called the perversity function. A locally closed set is an
intersection of an open and a closed set.
% If $C$ is a closed set and $U$ is
%an open set of $X$, then $ C\cap U$ is the difference of the two
%closed sets $C \cap (X - U) \subset C$.
%So partitions of a space into finitely many locally closed sets can
%be described by an increasing sequence of closed subsets.
The t-structure associated to a perversity function $p$ is given by
\[ \dd_{\geq 0} = \{ X \, | \, H_n (
i_{S_a}^* X ) =0 \text{ for } n < p (a ) , a \in A \} , \] for the locally
closed sets $S_a$ in $S$ \cite[2.2.1]{bbd}.
%If $S$ is a
%regular topological space (one can separate points and closed sets
%by open subsets), then t
This agrees with
\[ \{ X \, | \, H_n (
  X )_q =0 \text{ for } n < p (a ), q \in S_a , a \in A  \} . \]
 % where $Z ( S_a )$ is the interior of the intersection of all open
 % subsets of $S_{\text{sob}}$ containing $S_a$.
  Given a perversity function $p$. There is associated a corresponding
    function $ d_p \colon
S \rarr \zz \cup \{ \pm \infty \}$, defined by letting  \[
d_{p}^{-1} ( [n , \infty] ) - d_{p}^{-1} ( [n+1 , \infty] ) =
\cup_{a \in p^{-1} ( n) } S_a . \] The perverse t-structure associated
to $p$ agrees with the t-structure in Proposition \ref{t-structures}
for the associated function $d_p$. In particular, the perverse
t-structures on $\dd_R$ lift to t-model structures on $\ch ( \mm)$.

%We also get a t-model structure on the full subcategory of $\dd$
%consisting of bounded complexes.
% Perverse t-structures
%on categories of $R$-modules are often constructed using
% gluing techniques, rather than localization techniques.

%Assume that $S $ is a sober space and assume that $d \colon S \rarr
%\zz \cup\{ \pm \infty \}$ is a function such that
%is semicontinuous with respect to the ordering
%topology on $\zz \cup\{ \pm \infty \}$, i.e.~
%for each $n \in \zz \cup\{ \pm \infty \}$ the inverse image $ d^{-1}
%( \langle - \infty , n ] )$ is a closed subset of $ S$. Then
%$\dd_{\leq 0}$ equals \[ \{ X\, | \, H_n (X)_p = 0 \text{ whenever }
%n> d (p) \} . \] This follows since $R_U [m] \in \dd_{\geq 0}$ for
%all open subset $U$ of $S$ such that $U$ in $S$ is contained in the
%complement of $ d^{-1} ( m)$.
% Moreover, all complexes in $\dd_{\geq 0}$ can be built from
%the class of complexes of the form $R_U [m]$ described above.
\end{exmp}

\begin{exmp}[Flat model structure] \label{exmp:hov}
Let $ (S, R) $ be a ringed space, and assume
 that $ (S , R)$ has finite
hereditary global dimension \cite[3.1]{hov01}. Under these
assumptions Mark Hovey has constructed a
 tensorial model structure on $\ch ( \mm )$, called
  the flat model structure
   \cite[3.2]{hov01}.
%This is a technical
%assumption which implies that homology-isomorphisms of sheaves are in
%fact stalkwise homology-isomorphisms.
% (If all chain complex are quasi-isomorphic as presheaves to a
%complex of flasque sheaves.)
The assumptions are satisfied  if $S$ is a finite-dimensional
Noetherian space \cite[3.3]{hov01}. Hovey constructs a cofibrantly
generated model structure on $\ch (\mm )$ with weak equivalences the
stalkwise homology-isomorphisms. The fibrations are maps that are
levelwise surjective, as presheaves, and whose kernels are complexes
of flasque sheaves \cite[3.2]{hov01}.

We compare our model structure to his. Assume that $S =
S_{\text{sob}}$ and that
 $ (S , R)$ has finite
hereditary global dimension. Let $\ch (\nn) $ have the stalkwise model structure and $\ch (\mm)$ the flat model structure.
We claim that $L^2 \colon  \ch (\nn) \rarr \ch (\mm)$ is a Quillen left adjoint to the forgetfull functor $ i \colon \ch (\mm ) \rarr \ch (\nn)$. The functor $L^2$ respects cofibrations  since $L^2$ respects colimits and  $L^2$ of the cofibrant generators of the stalkwise model structure  are among
the   cofibrant generators of the flat model structure
  \cite[1.1, 3.2]{hov01}. In addition     $L^2$   respects weak equivalences. So $(L^2 , j)$ is a Quillen adjunction.
This is a  Quillen equivalence since $ L^2 X \to Y $ is a sheaf homology-isomorphism  in $\ch (\mm )$ if and only if the adjoint map $ X \to i (Y)$ is a sheaf homology-isomorphism in $\ch (\nn)$, for $ X \in \ch ( \nn ) $ and $ Y \in \ch ( \mm )$.
\end{exmp}

\begin{exmp}[Injective model structure]
  The injective model structure on $\ch (\mm)$ (or $\ch (\nn )$) has sheaf
homology isomorphisms as weak equivalences and levelwise injections
as cofibrations. For a discussion of this model structure see
 \cite[2.3.13]{hov}. The fibrant objects are chain complexes of
  injective
 sheaves
that are $K$-injective in the sense of Spaltenstein
\cite[1.1]{spal}.  Assume that $\ee$ has enough points.
 Then the  adjoint pair $(L^2, i)$ gives a Quillen equivalence between $ \ch (\nn)$ with the stalkwise model structure and $\ch (\mm)$ with the injective model structure. This follows as in the previous example because  $ L^2$ respects injections of presheaves.
  In particular, the class of
 fibrations for the stalkwise model structure on $\ch (\nn)$
    contains the  class of
 fibrations for the injective model structure on $\ch (\mm)$ (and $\ch (\nn)$).
%\begin{cor} There is a homological type spectral sequence with \[
%E^{2}_{s, t} = \textstyle\bigoplus_{k - l = t } Ext^s ( \hh_k (X) ,
%\hh_l (Y) ) \Rightarrow \dd ( X , Y )_{s + t } \] that converges
%weakly if the homology of $X$ is bounded below. The differential
%$d^r$ has degree $( -r , r -1 )$.
%\end{cor}

%\begin{lem} \label{lem:hohomstalkwise} Assume that $\ee$ is has
% enough points.
%If $X \in \ch
% (\nn)$ such that $ \hh_k (X) =0$ for $k < n$, then
% $ \dd ( R_C [n] , X ) $ is isomorphic to $
% L^2 \hh_n (X) (C)$. \end{lem}
%\begin{proof}
%\end{proof}
 \end{exmp}

\begin{exmp}[Stable homotopy category] We give an application of
Proposition \ref{t-structures} to a category which is not a derived
category. The heart of the stable homotopy category $\dd$ with the
 t-structure given by Postnikov sections is equivalent to
 the category of abelian
groups. Hence Proposition \ref{t-structures} gives a twisted variant
of the Postnikov t-structure. For each rational prime $p$ let $N_p
\in \zz \cup \{ \pm \infty \}$ and let $N_0 \in \zz \cup \{ \pm
\infty \}$ be greater or equal to $N_p$ for all primes $p$. Then the
associated full subcategory, $\dd_{\geq 0}$ (the connective
spectra), of $\dd$ consists of spectra $X$ such that
\[ (H_{n_p} (X))_p =0 \] whenever $n_p < N_p$, for all primes $p$,
and $H_{n_0} (X) \otimes \mathbb{Q} =0 $ whenever $ n_0 < N_0$.
%[Applied to the $G$-equivariant stable model category of orthogonal
%$G$-spectra we recover the t-structures constructed in \cite{gsp}.
%Note that the abelian category of coefficient systems for some
%$G$-universe $\mathcal{U}$ is a topos.]
%The coefficient systems are functors
%from the opposite stable $\uu$-category of orbit spectra to the
%category of abelian groups.
% The points are $\Phi G$.]
\end{exmp}

\begin{exmp}[Quasicoherent sheaves] \label{exmp:quasicoherent}
Let $(S, \oo_S)$ be a scheme.
  There is a full abelian
subcategory of $\mm$ consisting of quasi-coherent $\oo_S$-modules.   Typically, $\oo_U$
is not quasi-coherent for an open subset $U$ of $S$. So we can not follow Chapter \ref{sec:proj} and give a projective model structure on the category of chain complexes of
quasi-coherent $\oo_S$-modules.

We can construct  t-model structures on the category of chain complexes of quasi-coherent presheaves using  the techniques of Chapter \ref{sec:t-structures} and a proper cofibrantly generated model structure given by Mark Hovey  for certain  schemes \cite[2.4, 2.5]{hov01}.
The proof of Lemma \ref{simplicial} shows that this model structure
is quasi-simplicial (the complex associated to a simplicial set are
chain complexes of free $\oo_S$-modules).

With some assumptions on $S$ we can  inherit a t-structure on the
derived category of quasi-coherent $\oo_S$-modules from a t-structure on $\dd_{\oo_S}$.
 Assume that $S$ is a finite dimensional Noetherian scheme.
 Since $S$ is quasi-compact and quasi-separated   the
 derived category of chain complexes of quasi-coherent
 $\oo_S$-modules is a full
 subcategory of the derived category of chain
 complexes of $\oo_S$-modules \cite[p.187]{SGA6}.
 Moreover, our assumptions guarantee that
  the objects of this full
 subcategory are exactly complexes with quasi-coherent homology
  \cite[p.191]{SGA6}.
Hence the t-structures on $\dd_{\oo_S}$ constructed in Proposition
\ref{t-structuresrefined} restrict to give t-model structures on
the derived category of chain complexes of quasi-coherent
$\oo_S$-modules. See also the preprint by Roman
Bezrukavnikov \cite{bez}.

%These perverse t-structures come from t-model structures.
%T-structures on the derived category of chain complexes of
%quasi-coherent sheaves has been studied by
\end{exmp}

\begin{exmp}[Pro-chain complexes] A t-model structure on $
\ch (\nn)$ is  well suited to  give a  model structure on  the
category of pro-chain complexes of presheaves of $R$-modules. For
example there is a proper stable tensor model structure on $\p \ch
(\nn )$ such that the levelwise t-structure on the triangulated
homotopy category of $\p \ch (\nn )$ has the property that the
intersection
\[ \textstyle\bigcap_{n \in \zz } (\text{Ho} \, (\p \ch (\nn )))_{ \geq
n}
\]   consists of objects isomorphic to the
$0$-object \cite{ffi}. Recent work by the author on a general local
cohomology theory makes use of model structures on categories of
pro-chain complexes of presheaves of $R$-modules \cite{floc}.
\end{exmp}

%\begin{rem} One can construct a model structure on $\p \ch (\nn)$
% that is cofibrantly generated with concrete cofibrant generators,
%and such that the model structure restricts to give the stalkwise
%model structure on $ \ch (\nn)$. The key idea is to make use of the
%adjunction \[ \p \ch (\nn ) ( \{ R_C \} , N ) \cong \hh^0 ( N_x ) .
%\] [Check] The model structure is inherited from the right adjoint
%functor \[ \prod_{p \in \pt } \lim stalk_p .\]
%\end{rem}

\appendix
\section{Bousfield's cardinality argument} \label{appendix:bousfield}

 Bousfield's cardinality argument is used to
localize model categories \cite{bou75}. We give an extension of this
result. Let $\g$ be a cofibrantly generated model category. Let $I$
be a set of cofibrant generators. We assume that the sources of the
maps in $I$ are small.

\begin{defn} Let $X$ be an $I$-cell complex. The cardinality of the
set of cells in $X$ is denoted $\sharp X$. Let $i \colon A \rarr X$
be a relative $I$-cell complex. The cardinality of the set of relative
$I$-cells in the relative cell complex $i$ is denoted $ \sharp ( X ,
A )$.
\end{defn}

\begin{defn} Let $h$ be a functor from $\g$ to the category of sets.
%We assume that $h $ factors through the category of nonempty sets.
 The
functor $h$ is said to satisfy the \mdfn{colimit axiom} if for all
relative $I$-cell complexes $A \rarr X$
\[ \colim_{\alpha} \, h ( X_{\alpha} ) \rarr h (X ) \]
is a bijection where the colimit is over all relative sub $I$-cell complexes
$i_{\alpha} \colon A \rarr X_{\alpha}$ of $i$ such that $\sharp ( X_{\alpha }
, A) $ is finite.
\end{defn}
 The assumption that $\sharp ( X_{\alpha } , A) $
is finite can be relaxed. We restrict to the less general version as
that suffices for our examples.
\begin{defn} \label{defn:sigmauniform}
Let $\sigma$ be a cardinal number. We say that a functor $h$ from
  relative $I$-cell complexes to sets is \mdfn{$\sigma $-uniform}
if the cardinality of $ h (X)$ is less or equal to $ \sigma \times
\sharp X $ for all $I$-cell complexes $X$. \end{defn}

Given two sets of functors \mdfn{$ \{ i_a \}$} and \mdfn{$\{ s_b \} $} from $\g$ to the
category of sets, we define a class of maps in $\g$ depending on {$ \{ i_a \}$} and {$\{ s_b \} $}.
\begin{defn} Let \mdfn{$\ww$}
 denote the class of all maps $f \colon X
\rarr Y$ such that $i_a (f) $ is injective and $s_b (f) $ is surjective, for all $a$ and $b$.
\end{defn}
Note that the class $\ww$ is closed under composition and retract
 but that it need not satisfy the two-out-of-three property.
%We are
%interested in examples where $\ww$ does not satisfy the
%two-out-of-three property.
A typical example of a class $E$ is the class of $n$-equivalences associated to a
t-structure on the homotopy category of $\g$ (when $\g$ is a stable
model category).

%We could have included a third function $b$ that we required to be
%bijective, but it is more natural to consider $i \coprod b $ and $ s
%\coprod b$ instead.

We want to produce a (functorial) factorization of a map in $\g$ as
a map in $C \cap \ww$ followed by a map that has the right lifting
property with respect to all maps in $C \cap \ww$.

We say that $f \colon A \rarr X$ is an \mdfn{$I$-cell complex pair}
if $A$ is an $I$-cell complex and $f$ is a relative $I$-cell
complex.

\begin{defn} Let \mdfn{$ \ww_{\sigma} $}
 denote the class of all $I$-cell complex pairs $ X
\rarr Y$ with $\sharp Y\leq \sigma $ such that $i_a ( X ) \rarr i_a(Y)$
is injective and $ s_b (X ) \rarr s_b (Y)$ is surjective, for all $a$ and $b$. Let
\mdfn{$\ww'$} denote the class of all $I$-cell complex pairs $ X
\rarr Y$ such that $i_a ( X ) \rarr i_a (Y)$ is injective and $ s_b (X )
\rarr s_b (Y)$ is surjective, for all $a$ and $b$. Let \mdfn{$\ww''$} denote the class of
all relative $I$-cell complexes $ X \rarr Y$ such that $i_a ( X )
\rarr i_a (Y)$ is injective and $ s_b (X ) \rarr s_b (Y)$ is surjective, for all $a$ and $b$.
 \end{defn}

The class $ \ww_{\sigma} $ is skeletally small for each cardinal
$\sigma$, but $\ww'$ and $\ww''$ need not be skeletally small. In
the next Lemma one ought to be  careful about the meaning of
intersection sub-cell complexes. We follow Hirschhorn and assume that the model category is  cellular \cite[12]{hir}.

\begin{lem}
\label{bouskey} Let $\g$ be a cellular model category. Let $\sigma $
be an infinite cardinal number, and let $\{ i_a \}_{a \in A} $ and $\{ s_b \}_{b \in B}$ be
 $\sigma$-uniform functors from $\g$ to the category of sets. Assume that $A$ and $B$ have cardinality not greater than $\sigma$.

Let $ f \colon A \rarr X$ be an $I$-cell complex pair with $A \not=
X$. Assume that $f $ is in $ \ww'$. Then there is an $I$-cell
subcomplex $B$ of $X$ such that $ \sharp B \leq \sigma$, $B
\not\subset A$, and $ B\cap A \rarr B$ is in $\ww'$.
\end{lem}

\begin{proof}
 We construct an increasing
 sequence \[ B_0 \subset B_1 \subset B_2
 \subset \cdots \] of $I$-cell subcomplexes of $X$
 %indexed on an indexing set of cardinality $\sigma$
 such that:
 \begin{itemize}
 \item
$B_0 \not\subset A$.
\item whenever $i_1 , i_2 \in i_a (B_n \cap A ) $ map to the
same element  in $ i_a (B_n)$, then they map to the same element in $i_a
(B_{n+1} \cap A) $
\item the set $s_b ( B_n)$ maps to the image of $s_b ( B_{n+1} \cap A )$
in $ s_b ( B_{n+1 })$.
\end{itemize}

We choose some finite subcomplex $B_0$ of $X$ that is not contained
in $A$. We can do this since $ A \not=X$ and the gluing map from any
$I$-cell to $A$ factors through a finite  $I$-cell subcomplex of $A $ since
the sources of the maps in $I$ are small.

Assume that $B_n$ has been constructed. We construct $B_{n+1}$. Let
 $i_1 , i_2 \in i_a (B_n \cap A ) $ be two elements that map to the same
  element  in $i_a
(B_n)$. By our assumption on $f$ the two elements are sent to the same
element  under $i_a (B_n \cap A ) \rarr i_a ( A )$. The colimit axiom for
$i_a$ implies that there is a finite relative $I$-cell complex $B_n
\rarr I_{x_1 , x_2} $ in $X$ such that the two elements map to the
same element  in $ i_a ( I_{x_1 , x_2}\cap A ) $.

Similarly, for every element  $y \in s_b ( B_n)$ there is a finite
relative $I$-cell complex $ B_n \rarr S_y$ in $X$ such that the
image of $y$ in $s_b ( S_y )$ is in the image of $ s_b ( S_y \cap A)$.
Now let \[ B_{n+1} = B_n \cup_{i_1 , i_2} I_{i_1 , i_2 } \cup_y S_y
\] where the sum is over $i_1 , i_2 \in i (B_n \cap A ) $ and $y \in
s_b ( B_n)$, for all $a$ and $b$. This complex satisfies the conditions in the list above.

Now let $B $ be the union of all the $B_n$. The colimit axiom gives
that $ B \cap A \rarr B$ is in $\ww_{\sigma}$. The assumption that $i_a$, $s_b$ are $\sigma$-uniform, and that the cardinality of $A$ and $B$ are not greater than $\sigma$ give  that $\sharp B \leq \sigma$.
\end{proof}

\begin{lem} The
 class $\inj \ww' $ is equal to the class $
\inj \ww_{\sigma}$. \end{lem}
\begin{proof}
Let $ f \colon X \rarr Y$ be a map in $\ww'$.
 The trick is to use Lemma \ref{bouskey}
to write $ f $ as a transfinite composition of maps in $
\ww_{\sigma}$. There is a transfinite sequence $ X_{\lambda}$ such
that:
\begin{itemize}
\item
$ X_{\lambda} \rarr X_{\lambda +1 } $ is in $\ww_{\sigma}$
\item if
$\lambda$ is a limit ordinal, then $X_{\lambda} = \cup_{l < \lambda}
X_l $ \item if $X_{\lambda} $ is strictly contained in $X$, then
$X_{\lambda +1 } $ is strictly larger than $X_{\lambda} $.
\end{itemize}
Lemma \ref{bouskey} implies that we must have that $ X_{\lambda} =
X$ for some $\lambda$.
\end{proof}

\begin{lem} If $\ww''$ is closed under pushout along
cofibrations in $\g$, then
 $\inj \ww' = \inj \ww''$. \end{lem}
\begin{proof} This is a consequence of \cite[Prop.~13.2.1]{hir}.
\end{proof}

\begin{pro} \label{pro:bousfield}
There is a functorial factorization of the maps in $ \g $ with
an $I$-cell complex source as a map in $ \ww'$ followed by a map in
$ \inj \ww' $. Moreover, if $\ww'$ is closed under pushout along
cofibrations in $\g$, then there is a functorial factorization of
any map in $ \g $ as a map in $ \ww''$ followed by a map in $ \inj
\ww'' $.
\end{pro}
%\begin{rem}
%Proposition \ref{pro:bousfield} allow us to localize $ \text{Ho}
%(\g)$ with respect to many different functors. For example $ h_*
%\otimes_R M$ for any homology theory $h_*$ which satisfies the
%colimit axiom, and any module $M$. The module need not be flat so $
%h_* \otimes_R M$ need not be a homology theory.
%\end{rem}

\end{document}